\documentclass[12pt, a4paper]{article}
\usepackage[pdftex]{graphicx}

\newtheorem{theorem}{Theorem}

\begin{document}
\input{amssym}

\begin{center}
{\LARGE \bf Fundamental Dominations in Graphs}
\end{center}
\bigskip
\begin{center}
{\bf Arash Behzad\\ University of California,
LosAngeles\\abehzad@ee.ucla.edu}\\
\vspace{1.2cm} {\bf Mehdi Behzad\\Shahid Beheshti University,
Iran\\mbehzad@sharif.edu}\\
\vspace{1.2cm} {\bf Cheryl E. Praeger\\University of Western
Australia, Australia\\praeger@maths.uwa.edu.au}
\end{center}

\begin{center}
\bigskip\noindent
{\bf Abstract}
\end{center}
Nine variations of the concept of domination in a simple graph
are identified as fundamental domination concepts, and
a unified approach is introduced for studying them.
For each variation, the minimum cardinality
of a subset of dominating elements
is the corresponding fundamental domination number. It is observed that,
for each nontrivial connected graph, at most
five of these nine numbers can be different, and inequalities between
these five numbers are given. Finally, these fundamental
dominations are interpreted in terms of the total graph of the given graph, a concept
introduced by the second author in 1965. It is argued that the very first domination
concept, defined by O. Ore in 1962 and under a different name by
C. Berge  in 1958, deserves to be called the most fundamental
of graph dominations.\\


\bigskip
\section{Introduction and Preliminaries}

The literature contains extensive studies of many variations of the
concept of domination in a simple graph.
The following points are some loose extracts from the Preface,
Chapter 12, and the Appendix of the reference text [9], written
by well-known authorities: T.W. Haynes, S.T. Hedetniemi, and P.J.
Slater. \\

$\bullet$ One of the authors' objectives is ``{\it to consolidate
and organize much of the material in the more than 1200 papers
already published on domination in graphs}."

$\bullet$ ``{\it It is well known and generally accepted that the
problem of determining the domination number of an arbitrary
graph is a difficult one. ...  Because of this, researchers have
turned their attention to the study of classes of graphs for
which the domination problem can be solved in polynomial time}".

$\bullet$ ``{\it The following pages contain a fairly
comprehensive census of more than 75 models of dominating and
related types of sets in graphs which have appeared in the
research literature over the past 20 years.}"\\

\indent Without a doubt, the literature on this subject is growing
rapidly, and a considerable amount of work has been dedicated to
find different bounds for the domination numbers of graphs [3, 8,
10, 11]. The terms ``dominating set", and "domination number" of a
graph $G=(V,E)$ were first defined by O. Ore in 1962, see [12]. A
subset $A\subseteq V$ is a {\it dominating set} for $G$ if each
element of $V$ is either in $A$, or is adjacent to an element of
$A$. The {\it domination number} $\gamma(G)$, which is the most
commonly used domination number, is the minimum cardinality among
all dominating sets of $G$. We will interpret $\gamma(G)$ as the
minimum cardinality among all subsets $A\subseteq V$ dominating
the set $B=V$. Such a subset $A$ and the set $B$ are called a {\it
dominating set} and {\it the dominated set} for this {\it
vertex-vertex domination variation}, respectively. The parameter
$\gamma(G)$ will be referred to henceforth as the {\it
vertex-vertex domination number} of $G$.

\indent  Later, a few researchers defined, sometimes redefined,
and studied other domination variations: vertex-edge,
edge-vertex, etc., and the cardinalities of their largest or
smallest dominating sets [1,10].

\indent However, from practical point of view, it was necessary to
define other types of dominations. Most of these new variations
required the dominating set to have additional properties such
as: being as independent set, inducing a connected subgraph, or
inducing a clique. These properties were reflected in their names
as an adjective: independent domination, connected domination, and
clique domination, respectively.

\indent In this paper, among over 75 models of domination and
corresponding subsets in graphs, we choose nine variations, which
are the core, and call them {\it fundamental dominations}. For
this purpose, we consider simple nontrivial connected graphs
$G=(V,E)$. Our reason for restricting attention to such graphs
are two-fold. First, two of the fundamental domination numbers
are not defined for graphs with isolated vertices. Second, if $G$
is disconnected and has no isolated vertices, then the value of
each of the fundamental domination numbers for $G$ is equal to
the sum of the values of the same number for each of the
connected components of $G$.

\indent We introduce a unified approach to studying these fundamental
dominations  based on the fact that, in each domination variation,
two sets are used - the set consisting of the dominating elements,
and the set consisting of the elements that need to be dominated.
For each fundamental domination the
minimum cardinality among all dominating sets is the corresponding
{\it fundamental domination number}. We
observe, in Theorem 1, that for each nontrivial connected graph at most
five of these nine numbers can be different.
Inequalities concerning each pair of these five
numbers are considered in Theorems 2 and 3.\\

\indent Finally, we show how  these fundamental dominations may
be interpreted in terms of the total graph $T(G)$ of $G$,
defined by the second author in 1965. We argue that the very first domination
concept, defined by O. Ore in 1962 and under a different name by
C. Berge  in 1958, deserves to be called the {\it most fundamental
of graph dominations}.\\

\subsection{The fundamental domination numbers}

\indent By an {\it element} of a graph $G=(V,E)$ we mean a member
of the set $V\cup E$. Two different elements of $G$ are said to
be {\it associated} if they are adjacent or incident in $G$. For
$U, W\in\{V,E, V\cup E\}$, a subset $A\subseteq U$ dominates $W$
if each element of $W\backslash(A\cap W)$is associated with an
element of $A$. The minimum cardinality of such subsets $A$ is
denoted by $\gamma_{_{U,W}}(G)$. For historical reasons, we
replace $\gamma_{_{V,V}}(G)$, $\gamma_{_{E,E}}(G)$, and
$\gamma_{_{V\cup E,V\cup E}}(G)$ by $\gamma(G), \gamma^{'}(G)$,
and $\gamma^{''}(G)$, respectively. We call the following nine
parameters the {\it fundamental domination numbers} of $G$:
$$\gamma_{_{V,V}}(G)=\gamma(G); \ \ \gamma_{_{V,E}}(G); \ \ \gamma_{_{V,V\cup E}}(G);$$
$$\gamma_{_{E,V}}(G); \ \ \ \ \gamma_{_{E,E}}(G)=\gamma^{'}(G); \ \ \gamma_{_{E,V\cup E}}(G);$$
$$\gamma_{_{V\cup E,V}}(G); \gamma_{_{V\cup E,E}}(G); \ \ \ \gamma_{_{V\cup E,V\cup E}}(G)=\gamma^{''}(G).$$
In each case, a dominating subset whose cardinality equals the
fundamental domination number is called a {\it fundamental
dominating subset} of $G$.

\indent Various equalities and inequalities related to these
fundamental numbers are summarized in Theorems 1, 2, and 3 of
Section 2. See Figure 2 which actually reflects the statements of
these three theorems.

\indent The {\it line graph} $L(G)$ of a nonempty graph $G=(V,E)$
is the graph whose vertex set is in one-to-one correspondence with
the elements of the set $E$ such that two vertices of $L(G)$ are
adjacent if and only if they correspond to two adjacent edge of
$G$. The {\it total graph} $T(G)$ of $G$ is the graph whose vertex
set is in one-to-one correspondence with the set $V\cup E$ of
elements of $G$ such that two vertices of $T(G)$ are adjacent if
and only if they correspond to two adjacent or incident elements
of $G$, see [4,6]. Figure 1 shows a graph $G$, its line graph
$L(G)$ and its total graph $T(G)$. It is clear that $G$ and $L(G)$
are disjoint induced subgraphs of $T(G)$.\\


\begin{figure}
\begin{center}
\includegraphics[width=4.5in, clip=false, trim=13 4.5in 10 5in]{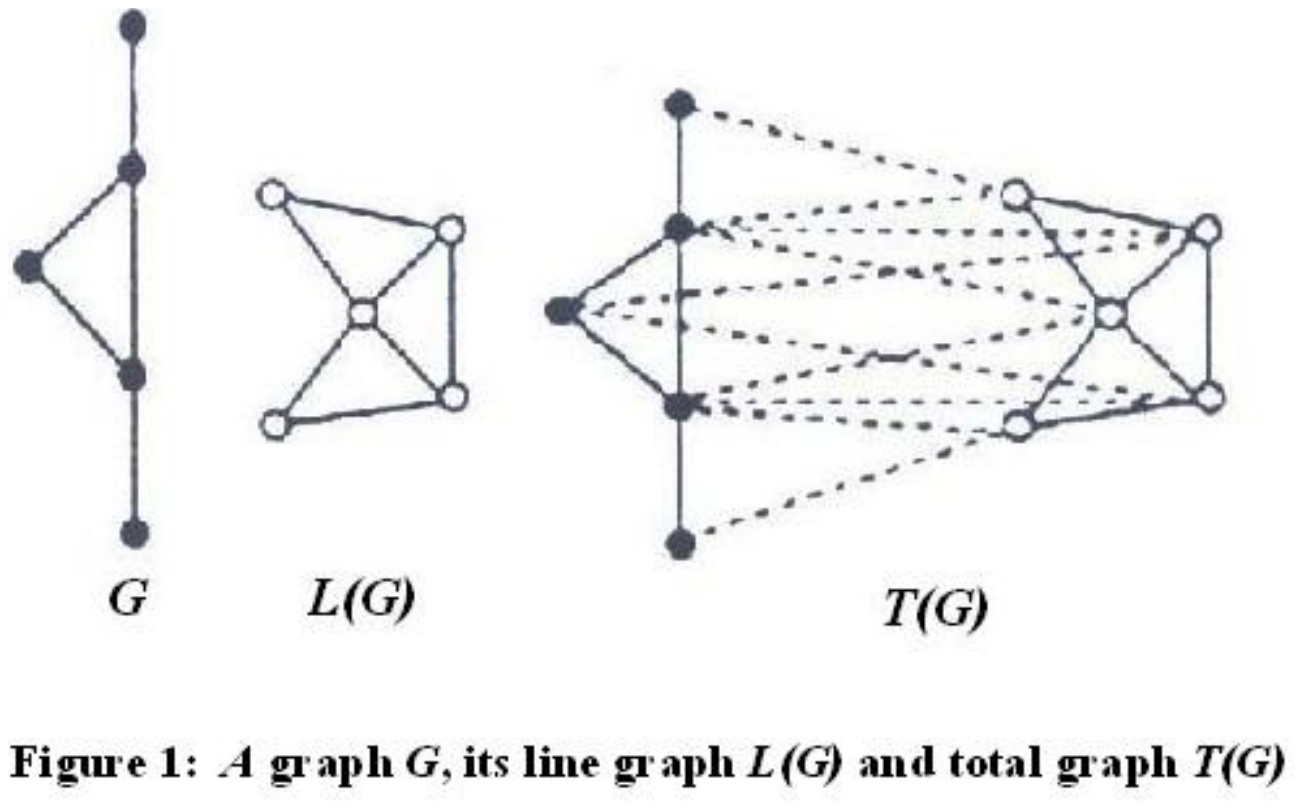}
\end{center}
\end{figure}

\indent In Section 3, we elaborate on the fact that each graph
domination variation can be presented in the context of a special
subset $A$ of $V(T(G))$ dominating an appropriate subset $B$ of
$V(T(G))$. This observation suggests that the phrase ``the most
fundamental of graph dominations" might appropriately be attached
to the vertex-vertex domination related to $\gamma$.

\bigskip
\indent Notions and notations not defined here can be found in
texts such as M. Behzad, et al [6], and D. B. West [14].

\bigskip
\section{General Equalities and Inequalities}
\indent Each of the nine fundamental domination numbers has its
own applications, and must be considered separately. However, for
each graph under consideration the following equalities hold.

\begin{theorem}
{\it For each nontrivial connected graph} $G$ {\it we have}:

 \noindent (1) $\gamma^{'}(G)=\gamma(L(G)),$ \ \
\ \ (2) $\gamma^{''}(G)=\gamma(T(G)),$ \ \ \ \ \ \ \ (3)
$\gamma_{_{V\cup E,V}}(G)=\gamma(G)$,\\
(4) $\gamma_{_{V\cup E,E}}(G)=\gamma^{'}(G)$, \ \ (5)
$\gamma_{_{E,V}}(G)=\gamma_{_{E,V\cup E}}(G)$, \ \ (6)
$\gamma_{_{V,E}}(G)=\gamma_{_{V,V\cup E}}(G).$
\end{theorem}

\noindent {\bf Proof.} The equalities (1), and (2) follow from the
definitions of the line graph and the total graph of $G$,
respectively. To prove the last four equalities, first we observe
that, by definition, the parameter on the left of each equation
in less than or equal to the parameter on the right.\\
\indent For example, in (3) we show that $\gamma_{_{V\cup
E,V}}(G)\leq\gamma(G)=\gamma_{_{V,V}}(G)$. Since each subset of
$V$ dominating $V$ is a subset of $V\cup E$ that dominates $V$,
the inequality $\gamma_{_{V\cup
E,V}}(G)\leq\gamma_{_{V,V}}(G)=\gamma(G)$ follows. The same is
true for the inequality $\gamma_{_{V\cup
E,E}}(G)\leq\gamma^{'}(G)=\gamma_{_{E,E}}(G)$.

\indent Next, we prove that
$\gamma_{_{E,V}}(G)\leq\gamma_{_{E,V\cup E}}(G)$. Subsets such as
$A$ of the set $E$ that dominates $V$ do not necessarily
dominates $V\cup E$. Ordinarily, subsets of bigger size are
needed to dominate both $V$ and $E$. Thus by definitions involved
$\gamma_{_{E,V}}(G)\leq\gamma_{_{E,V\cup E}}(G)$. The same
argument shows that $\gamma_{_{V,E}}(G)\leq\gamma_{_{V,V\cup E
}}(G)$.\\
\indent To complete the proofs of the equalities (3)-(6) it
suffices to prove the validity of each converse inequality. As an
example, for (3) we show that $\gamma(G)\leq\gamma_{_{V\cup
E,V}}(G)$. Let $A$ be a subset of $V\cup E$ which dominates $V$
such that $|A|=\gamma_{_{V\cup E,V}}(G)$. If $A\cap E=\phi$, then
$A\subseteq V$, and $\gamma(G)\leq\gamma_{_{V\cup E,V}(G)}$.
Suppose $A\cap E\neq\phi$. For each $e=uv\in A\cap E$ eliminate
$e$ from $A$ and, if necessary, add to the remaining subset one of
$u$ or $v$ to produce a set $A^{'}\subseteq V$ dominating $V$.
Thus $\gamma(G)\leq|A^{'}|\leq|A|=\gamma_{_{V\cup E,V}}(G)$.

\indent Next, we prove that $\gamma^{'}(G)\leq\gamma_{_{V\cup
E,E}}(G)$. Let $A\subseteq V\cup E$ which dominates $E$ such that
$|A|=\gamma_{_{V\cup E,V}}(G)$. If $A\cap V=\phi$, then
$A\subseteq E$, and $\gamma^{'}(G)\leq\gamma_{_{V\cup E,E}}(G)$.
Otherwise, $A\cap V\neq\phi$. For each $v\in A\cap V$, eliminate
$v$ from $A$, and if necessary, add to the remaining subset an
edge incident with $v$ to produce a set $A^{'}\subseteq E$
dominating $E$. Such an edge exists, since $G$ has no isolated
vertices. Thus, as before,
$\gamma^{'}(G)\leq|A^{'}|\leq|A|=\gamma_{_{V\cup E,V}}(G)$, and
(4) is established.

\indent To Prove $\gamma_{_{E,V\cup E}}(G)\leq\gamma_{_{E,V}}(G)$,
let $A\subseteq E$ dominates $V$, such that
$|A|=\gamma_{_{E,V}}(G)$. We claim that $A$ dominates $E$ as well.
Let $e=uv\in E\backslash A$. Since $A$ dominates $V$, there
exists an edge $e^{'}\in A$ such that $e^{'}$ and $u$ are
incident, that is to say $u\in e^{'}$. Then $e$ and $e^{'}$ are
adjacent. Thus $A$ dominates $V\cup E$. Hence $\gamma_{_{E,V\cup
E}}(G)\leq\gamma_{_{E,V}}(G)$.

\indent Finally, in a similar manner, we prove that
$\gamma_{_{V,V\cup E}}(G)\leq\gamma_{_{V,E}}(G)$. Assume that
$A\subseteq V$ dominates $E$, and that $|A|=\gamma_{_{V,E}}(G)$.
We show that $A$ dominates $V$, too. Let $v\in V\backslash A$.
Since $G$ is connected and  nontrivial, there exists at least one
edge $e=uv$ incident with $v$, and since $A$ dominates $E$, one
of the vertices $u,v$ must lie in $A$. Since $v\notin A, u$ must
be in $A$. Hence $A$ dominates $V\cup E$. Therefore,
$\gamma_{_{V,V\cup
E}}(G)\leq|A|=\gamma_{_{V,E}}(G).\hspace{7.5cm}\square$\\

\indent Based on Theorem 1, the values of the nine fundamental
numbers are essentially reduced to five:
$$\gamma=\gamma_{_{V\cup E,V}}, \gamma_{_{V,E}}=\gamma_{_{V,V\cup E}}, \gamma_{_{E,V}}=\gamma_{_{E,V\cup E}}, \gamma^{'}=
\gamma_{_{V\cup E,E}}, \ and  \ \gamma^{''}.$$
\indent Next, we introduce a digraph $D$ with node set
 $N=\{\gamma,\gamma^{'},\gamma^{''},\gamma_{_{V,E}},\gamma_{_{E,V}}\}$,
 and arc set
\begin{eqnarray*}
A=\{(\gamma_{_{V,E}},\gamma^{'}),(\gamma_{_{V,E}},\gamma^{''}),(\gamma_{_{V,E}},\gamma),(\gamma_{_{V,E}},\gamma_{_{E,V}}),
(\gamma_{_{E,V}},\gamma_{_{V,E}}),\\
(\gamma^{'},\gamma),(\gamma,\gamma^{'}),(\gamma_{_{E,V}},\gamma^{'}),(\gamma^{''},\gamma^{'}),(\gamma_{_{E,V}},\gamma^{''})
,(\gamma^{''},\gamma),(\gamma_{_{E,V}},\gamma)\}.
\end{eqnarray*}
See Figure 2. In this figure the two arcs $(\gamma,\gamma^{'})$,
and $(\gamma^{'},\gamma)$ are represented by the straight line
joining the two nodes $\gamma$ and $\gamma^{'}$ along with two
arrows in opposite  directions. This situation will be denoted by
$\gamma\leftrightarrow\gamma^{'}$, or by
$\gamma^{'}\leftrightarrow\gamma$. The two arcs
$(\gamma_{_{V,E}},\gamma_{_{E,V}})$ and
$(\gamma_{_{E,V}},\gamma_{_{V,E}})$ are represented in the same
way. However, for each of the remaining eight arcs $(x,y)$ in $A$
such that $(y,x)\notin A$, the arc $(x,y)$ is represented by a
straight line joining the two nodes $x$ and $y$ along with one
arrow from $x$ to $y$. This situation will be denoted by
$x\rightarrow y$.


\begin{figure}
\begin{center}
\includegraphics[width=4in, clip=false, trim=7 4in 7 5in]{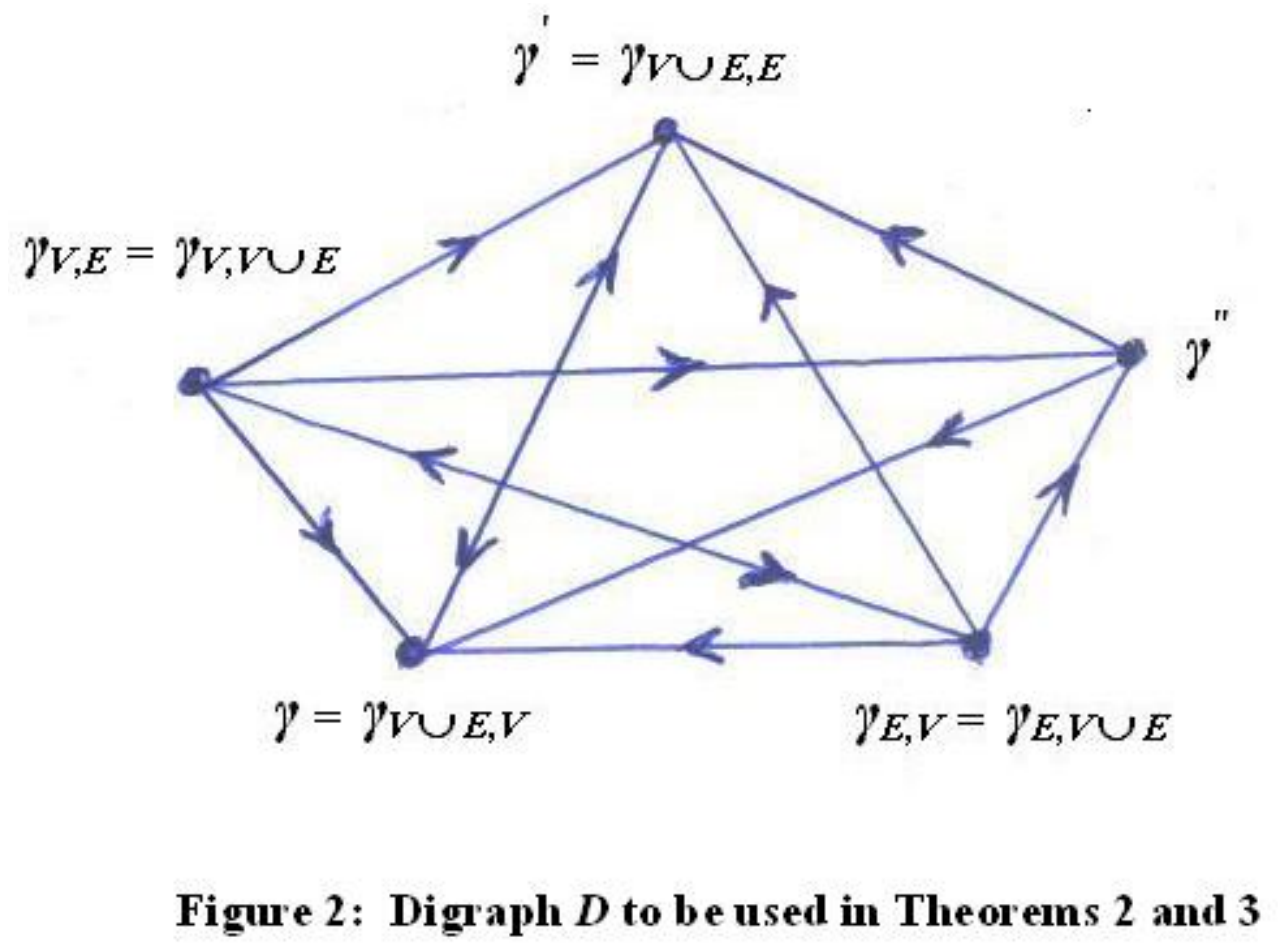}
\end{center}
\end{figure}

\begin{theorem}
Let $x$ and $y$ be two different nodes of the digraph $D$ of
Figure 2. If $(x,y)$ is an arc of $D$ and $(y,x)$ is not an arc
of $D$, then for every nontrivial connected graph $G$, the
inequality $x(G)\geq y(G)$ holds.
\end{theorem}
{\bf Proof.} The equalities stated in Theorem 1 are depicted in
Figure 2. In the light of these equalities, it suffices to show
that the following eight inequalities are valid:
$$
\begin{array}{lcclcc}
  &(1)& \gamma_{_{V,E}}(G)\geq\gamma_{_{V\cup E,E}}(G), ~~~~  &(2)& \gamma_{_{V,V\cup E}}(G)\geq\gamma^{''}(G), \\
&(3)& \gamma_{_{E,V}}(G) \geq\gamma_{_{V\cup E,V}}(G), ~~~~~&(4)& \gamma_{_{E,V\cup E}}(G)\geq\gamma^{''}(G), \\
  &(5)& \gamma_{_{V,V\cup E}}(G)\geq\gamma(G), ~~~~~~~ &(6)&
\gamma^{''}(G)\geq \gamma_{_{V\cup E,V}}(G), \\
&(7)& \gamma_{_{E,V\cup E}}(G)\geq\gamma^{'}(G), ~~~~~~~ &(8)&
\gamma^{''}(G)\geq\gamma_{_{V\cup E,E}}(G).
\end{array}
$$

\indent Proofs of the validity of each of the first four
inequalities are similar. ``Having more freedom to choose the
required dominating sets" is the key phrase. Note that in both
sides of each of these inequalities the same set needs to be
dominated; however, each dominating set of the parameter on the
left is necessarily a dominating set of the parameter on the
right too. As an example, (2) holds, since
$\gamma^{''}=\gamma_{_{V\cup E,V\cup E}},V\cup E=V\cup E$, and
$V\subseteq V\cup E$.

Similar arguments prove the last four inequalities. For each
inequality, the dominating sets of the two parameters which
appear on either side are subsets of the same set; however, for
the parameter on the right the set which needs to be dominated is
a subset of the same set for the parameter on the left. As an
example, (7) holds since $\gamma^{'}=\gamma_{_{E,E}}, E=E$, and
$E\subseteq V\cup E.\hspace{3.5cm}\square$\\

In the sequel, we will refer to the following special classes of
graphs, along with some of their specified domination numbers. We
will then use these considerations in Theorem 3 to explain the
significance of the two ``double arcs" $x\leftrightarrow y$ in
Figure 2.\\

\noindent{\bf Example 1.} For a complete graph, $K_{n}$, of order
$n, \ n\geq 2$, one can readily see that:
$$\gamma(K_{n})=1; \ \ \ \ \ \gamma^{'}(K_n)=\left\lfloor\frac{n}{2}\right\rfloor ;$$
$$\gamma_{_{E,V}}(K_{n})=\left\lceil {\frac{n}{2}} \right\rceil ; \ \
\ \ \gamma_{_{V,E}}(K_{n})=n-1.
$$

{\bf Example 2.} It is easy to observe that for the complete
bipartite graph $K_{1,n},n\in\Bbb N$, we have
$$\gamma_{_{V,E}}(K_{1,n})=1; \ \ \ \gamma_{_{E,V}}(K_{1,n})=n.$$

{\bf Example 3.} Next we introduce a new class of graphs,
$R_{3n},n\in\Bbb N$, called {\it ridged graphs}. Let:
$$V(R_{3n})=\{u_{i},v_{i},w_{i}|1\leq i\leq n\},  \ \ and$$
$$E(R_{3n})=E_{1}\cup E_{2}\cup E_{3}, \ \ where$$
$$E_{1}=\{u_{i}u_{i+1}|1\leq i\leq n-1\}, \ \ E_{2}=\{u_{i}v_{i}|1\leq i\leq n\}, \ \ E_{3}=\{u_{i}w_{i}|1\leq i\leq n\}.$$
The graph $R_{18}$ is shown in Figure 3.

\begin{figure}[tbh]
\begin{center}
\setlength{\unitlength}{0.4in}
\begin{picture}(5,2.5)

\put  (0,     0) {\makebox(0,0){$\bullet$}}

\put  (1,     0) {\makebox(0,0){$\bullet$}}

\put  (2,     0) {\makebox(0,0){$\bullet$}}

\put  (3,     0) {\makebox(0,0){$\bullet$}}

\put  (4,     0) {\makebox(0,0){$\bullet$}}

\put  (5,     0) {\makebox(0,0){$\bullet$}}

\put  (0,     1) {\makebox(0,0){$\bullet$}}

\put  (0,     2) {\makebox(0,0){$\bullet$}}

\put  (1,    2) {\makebox(0,0){$\bullet$}}

\put  (2,     2) {\makebox(0,0){$\bullet$}}

\put  (3,     2) {\makebox(0,0){$\bullet$}}

\put  (1,     1) {\makebox(0,0){$\bullet$}}

\put  (4,     2) {\makebox(0,0){$\bullet$}}

\put  (5,     2) {\makebox(0,0){$\bullet$}}

\put  (2,     1) {\makebox(0,0){$\bullet$}}

\put  (3,     1) {\makebox(0,0){$\bullet$}}

\put  (4,     1) {\makebox(0,0){$\bullet$}}

\put  (5,     1) {\makebox(0,0){$\bullet$}}

 \put  (-.105,    -.3)          {$w_1$}
\put  (.895,    -.3)          {$w_2$}

 \put  (1.895,-.3)             {$w_3$}

\put  (2.895,-.3)             {$w_4$}

\put  (3.895,-.3)             {$w_5$}

\put  (4.895,-.3)             {$w_6$}

\put  (-.105,    2.2)          {$v_1$}

\put  (.895,    2.2)          {$v_2$}

\put  (1.895,    2.2)          {$v_3$}

\put  (2.895,    2.2)          {$v_4$}

 \put  (3.895,    2.2)          {$v_5$}

 \put  (4.895,    2.2)          {$v_6$}

\put  (0.1,    1.1)          {$u_1$}

\put  (1.1,    1.1)          {$u_2$}

\put  (2.1,    1.1)          {$u_3$}

\put  (3.1,    1.1)          {$u_4$}

 \put  (4.1,    1.1)          {$u_5$}

 \put  (5.1,    1.1)          {$u_6$}

 \put(0,0)          {\line(0,  1) {2}}
 \put(1,0)          {\line(0,  1) {2}}
 \put(2,0)          {\line(0,  1) {2}}
 \put(3,0)          {\line(0,  1) {2}}
 \put(4,0)          {\line(0,  1) {2}}
 \put(5,0)          {\line(0,  1) {2}}
\put(0,1)          {\line(1,  0) {5}}

\end{picture}\\
Figure 3\quad The ridged graph $R_{18}$
\end{center}
\end{figure}

By an easy induction on $n$ one can prove that:
$$\gamma(R_{3n})=n; \ \ and  \ \gamma^{'}(R_{3n})=\left\lceil {\frac{n}{2}} \right\rceil.$$
\begin{theorem}
If $x\leftrightarrow y$ is a double arc of the digraph $D$ of
Figure 2, then for each positive integer $r$ there exist graphs
$G$ and $H$ such that $x(G)-y(G)>r$, and $y(H)-x(H)>r$.
\end{theorem}

\noindent {\bf Proof.} The digraph $D$ contains two double arcs
of the form $x\leftrightarrow y$: namely
$\gamma^{'}\leftrightarrow\gamma$, and
$\gamma_{_{V,E}}\leftrightarrow\gamma_{_{E,V}}$. Hence, for a
given $r$ we must provide four graphs $G_{1},G_{2},G_{3}$, and
$G_{4}$ such that
$$\gamma^{'}(G_{1})-\gamma(G_{1})>r;$$
$$\gamma(G_{2})-\gamma^{'}(G_2)>r;$$
$$\gamma_{_{V,E}}(G_{3})-\gamma_{_{E,V}}(G_{3})>r;$$
$$\gamma_{_{E,V}}(G_{4})-\gamma_{_{V,E}}(G_{4})>r.$$
\indent To present $G_{1}$, let $n=2r+4$, and $G_{1}=K_{n}$. As
specified in Example 1, we have:
$\gamma^{'}(K_{n})-\gamma(K_{n})=\left\lfloor\frac{2r+4}{2}\right\rfloor-1=r+1>r$.

For $G_{2}$, let $n=4r$, and $G_{2}=R_{3n}$. Then Example 3
indicates that: $\gamma(R_{3n})-\gamma^{'}(R_{3n})=n-\left\lceil
{\frac{n}{2}} \right\rceil=2r>r$.

To provide $G_{3}$, as for $G_{1}$, let $n=2r+4$, and
$G_{3}=K_{n}$. Then, by Example 1:
$\gamma_{_{V,E}}(K_{n})-\gamma_{_{E,V}}(K_{n})=n-1-\left\lceil
{\frac{n}{2}} \right\rceil=r+1>r$.

Finally let $n=r+2$, and $G_{4}=K_{1,n}$. Then, Example 2 implies
that:
$\gamma_{_{E,V}}(K_{1,n})-\gamma_{_{V,E}}(K_{1,n})=n-1=r+1>r\hspace{5cm}\square$

\section{The Most Fundamental of Graph Dominations}
The total graph $T(G)$ of a nonempty graph $G=(V,E)$ was defined
in Section 1. It has two vertex disjoint subgraphs $G^{\ast}$ and
$L^{\ast}$ such that $G^{\ast}$ and $L^{\ast}$ are, respectively,
isomorphic to $G$ and to the line graph $L(G)$. See Figure 1, and
for a characterization of total graphs see [5].

Each graph domination concept can be presented in the context of
a special subset $A$ of the vertex set $V\cup E$ of $T(G)$
dominating an appropriate subset $B$ of $V\cup E$. For the nine
Fundamental domination variations these special subsets $A$ and
$B$ are specified below:\\

\indent ~~~~$A\subseteq V(G^{\ast})\subset V(T(G))$  \ for \
$\gamma(G), \gamma_{_{V,E}}(G)$, \ \ \ and \ $\gamma_{_{V,V\cup E}}(G)$;\\
\indent ~~~~$A\subseteq V(L^{\ast})\subset V(T(G))$  \ for \
$\gamma_{_{E,V}}(G), \gamma^{'}(G)$, \ and \ $\gamma_{_{E,V\cup
E}}(G)$;\\
\indent ~~~~$A\subseteq V(T(G))$ \hspace{1.4cm}  for \
$\gamma_{_{V\cup E,V}}(G), \gamma_{_{V\cup E,E}}(G)$, \ and \
$\gamma^{''}(G)$;\\
\indent ~~~~$B=V(G^{\ast})\subset V(T(G))$  \ for \
$\gamma(G), \gamma_{_{E,V}}(G)$, \ \ \ and \ $\gamma_{_{V\cup E,V}}(G)$;\\
\indent ~~~~$B=V(L^{\ast})\subset V(T(G))$  \ for \
$\gamma_{_{V,E}}(G), \gamma^{'}(G)$, \ and \ $\gamma_{_{V\cup
E,E}}(G)$;\\
\indent ~~~~$B=V(T(G))$ \hspace{1.4cm}  for \ $\gamma_{_{V,V\cup
E}}(G), \gamma_{_{E,V\cup E}}(G)$, \ and \
$\gamma^{''}(G)$.\\

Hence, for any fundamental domination number of $G$ one can
simply minimize the cardinalities of special subsets of $V(T(G))$
that dominate an appropriate subset of $V(T(G))$. Thus a slight
modification of the original concept defined in 1962 by Ore [12],
and in 1958 under a different name by Berge [7], deserves to be
referred to as the most fundamental of graph dominations.

\section*{Concluding Remarks}
\indent 1. In this paper we have attempted to categorize
domination concepts into nine categories. For a nontrivial
connected graph $G=(V,E)$, if we use $V$ for the set of vertices,
and $E$ for the set of edges of $G$,then the nine categories
might be referred to as dominations in: $V-V$, $V-E$, $V-V\cup E,
E-V, E-E, E-V\cup E, V\cup E-V, V\cup E-E$, and $V\cup E-V\cup E$
contexts. The majority of the concepts and results in this vast
area are related to $V-V$ dominations; see the book by S. T.
Hedetniemi and R. C. Laskar [11]. Many of the concepts and
results can easily be transformed into other categories. As an
example, we refer to the $k-$domination defined in the $V-V$
context as follows. For a positive integer $k$, a subset
$S\subseteq V$ is a {\it k-dominating} set if for each $u\in V-S,
\ |N(u)\cap S|\geq k$, where $N(u)$ denotes the set of vertices
of $G$ which are adjacent to $u$. The $k-$domination number
$\gamma_{_{k}}(G)$ was considered by E. J. Cockaye, B. Gamble, and
B. Shepherd in [8]. Considering this parameter in other contexts
and obtaining their values seems to be of interest.

2. Theorems 1, 2, and 3 provide relationships among the nine
fundamental domination numbers of a nontrivial connected graph
$G$. One can often combine some of the given bounds with some of
the existing results to produce new bounds which might, or might
not be sharp. As an example we consider
$\gamma^{''}(G)=\gamma(T(G))$, and the following result of C.
Payan [13]: if $G$ has order $n$, size $m$, and minimum degree
$\delta$, then $\gamma(G)\leq(n+2-\delta)/2$. Since the graph
$T(G)$ has order $n+m$, and its minimum degree is $2\delta$, we
have $\gamma^{''}(G)\leq 1-\delta+(n+m)/2$. Hence Theorems 1, and
2 imply that for a nontrivial connected graph $G$ we have:
$$\gamma^{'}(G)=\gamma_{_{V\cup E,E}}(G)\leq\gamma^{''}(G)=\gamma(T(G))\leq 1-\delta+(n+m)/2.$$
\section*{Acknowledgements}
The first two authors are thankful for being able to participate
in the Kashkul Project led by Professor E. S. Mahmoodian [2]. All
the authors would like to thank Babak Behzad for his technical
assistance.

\end{document}